\magnification1200

\newcount\secno\secno=0
\newcount\thmno\thmno=0
\def\tag{\the\secno.\the\thmno}

\def\section #1\par{
\bigskip
\global\advance\secno by 1\global\thmno=0{\tenfrak\the\secno.} {\bf #1}\par
\smallskip}

\def\proclaimno #1. #2\par{\global\advance\thmno by 1%
{\bf #1 \tag.} {\it #2}\par}

\font\tener=eurm10
\font\sevener=eurm7
\font\fiveer=eurm5
\newfam\eufam
\textfont\eufam=\tener
\scriptfont\eufam=\sevener
\scriptscriptfont\eufam=\fiveer
\def\eu#1{{\fam\eufam #1}}

\font\tenfrak=eufm10
\font\sevenfrak=eufm7
\font\fivefrak=eufm5
\newfam\frakfam
\textfont\frakfam=\tenfrak
\scriptfont\frakfam=\sevenfrak
\scriptscriptfont\frakfam=\fivefrak
\def\frak{\fam\frakfam}
\def\frakN{{\frak N}}

\font\tenams=msbm10
\font\sevenams=msbm7
\font\fiveams=msbm5
\newfam\amsfam
\textfont\amsfam=\tenams
\scriptfont\amsfam=\sevenams
\scriptscriptfont\amsfam=\fiveams
\def\Bbb#1{{\fam\amsfam #1}}

\def\sdp{\Bbb n}
\font\amsfa=msam10
\def\qed{\hfill\hbox{\amsfa\char'003}}
\def\Cal#1{{\cal #1}}

\font\ninerm=cmr8
\font\ninebf=cmb8
\font\tt=cmtt10

\nopagenumbers
\font\sevenss=cmss7
\headline{\sevenss \ifnum\pageno>1
\ifodd\pageno
a-T-menability of groups acting on trees\hrulefill\ \the\pageno
\else\the\pageno\ \hrulefill\ S.~R.~Gal
\fi\fi}

\parskip=\smallskipamount
\parindent=0pt

\centerline{\bf a-T-menability of groups acting on trees}
\footnote{}{2000 {\it Mathematics Subject Classification :} 20F65}
\centerline{\'Swiatos\l aw R. Gal\footnote{$^\star$}{Partially supported
by a KBN grant 2 P03A 017 25.}}
\centerline{Wroc\l aw University}
\centerline{\tt http://www.math.uni.wroc.pl/\~{}sgal/}
\bigskip

{\advance\leftskip by \bigskipamount\advance\rightskip by \bigskipamount

{\ninebf Abstract.} \ninerm We present some partial results
concerning a-T-menability of groups acting on trees. Various known
results are given uniform proofs.

}
\bigskip

\section Definition of a-T-menability.

\proclaimno Definition. Given a metric space $(X,\rho)$, the action
$\eu G\to{\rm Isom}(X)$ is {\rm metrically proper} if given $x\in X$
the displacement function $\eu G\ni g\mapsto\rho(x,gx)\in\Bbb R$ is proper.
If the group action is set we say that $X$ is metrically proper.

The above property does not depend on a choice of a point $x$.

\proclaimno Definition (M.~Gromov). A locally compact,
second countable, compactly generated group $\eu G$ is {\rm a-T-menable}
if there exists
metrically proper isometric  $\eu G$-action on some affine Hilbert space. 

{\it Remark :} a-T-menability is equivalent to the existence of
$C_0$-approximate unity consisting of positive definite functions.
The latter is called in the literature the Approximation Property
of Haagerup.

Throughout this paper by representation we mean isometric affine action.
By a subgroup we always mean a closed subgroup.

\section Motivation.

We are motivated by the following Theorems:

\proclaimno Theorem. [JJV, Th.~6.2.8] Let $\Gamma$ be a countable group
acting on a tree without inversions, with finite edge stabilizers. If vertex
stabilizers in $\Gamma$ are a-T-menable, then so is $\Gamma$.

\edef\jjv{Theorem \tag\ }

\proclaimno Theorem. [JJV, Ex.~6.1.6] If $\eu N$ is a-T-menable
and $\eu G/\eu N$ is amenable then $\eu G$ is a-T-menable.

\edef\joli{Theorem \tag\ }

On the other hand:

\proclaimno Theorem. [HV, Ch.~8 L.~6] Every isometric affine representation of
$G=SL_2\Bbb Z\sdp\Bbb Z^2$ has a $\Bbb Z^2$-fixed point.

The group $\eu G$ acts on a tree as it
can be decomposed as
$(\Bbb Z/4\sdp\Bbb Z^2)*_{\Bbb Z/2\sdp\Bbb Z^2}(\Bbb Z/6\sdp\Bbb Z^2)$.
The factors are a-T-menable by \joli, since
$\Bbb Z^2$ is a-T-menable as it acts on $\Bbb R^2$.
On the other hand
$Sl_2\Bbb Z$ is a-T-menable by the result of Haagerup [H].
Therefore extra assumptions about finitness of edge stabilizers (in \jjv)
or amenability of the quotients (in \joli) cannot be simply weakened to
a-T-menability.

We believe that a-T-menability is a property
of representations (see next section for definitions),
rather than that of groups.
Therefore in this mood we will concentrate on a question: what are the
conditions under which given proper affine representations of two
groups extend to one of their product with amalgamation, rather than
whether there exist some such representation that extend.

The most na\"\i ve observation is that
if there is a proper representation of a free product with amalgamation,
the restrictions are proper representations of the factors that coincide
on the common subgroup. It is not known whether the reverse holds,
however we will prove some results in this direction.
We would like to thank Tadeusz Januszkiewicz for calling our attention to it.

The main result of this paper is Theorem 5.4.
As a result we strengthen a result from [JJV] (see Section 6).
Constructions, we give, are purely geometric.
We will also give an affirmative answer to the question
of A.~Valette, whether Baumslag-Solitar groups are a-T-menable.
This was originally done (using another approach) in [GJ].

Finally we would like to thank Agnieszka for her hospitality in Vienna
and Jan Dymara for carefully reading
preliminary version of this paper
and his assistance in improving the presentation.

\section Actions on trees in general.

Before we examine the case of a free product with amalgamation,
let us restate a general observation 
of Haagerup [H]
and its easy generalizations.

Let $T$ be a tree. By $E$ we denote a vector space of functions
on edges of T with finite support. By $V$ we denote the affine space
of functions on vertices of $T$ with finite support and total mass one.
The structure of affine space is given as follows: $\delta_v-\delta_w$
is equal to the characteristic function of the segment joining vertices
$v$ and $w$ (with the appropriate signs with respect so some auxiliary
orientation on the edges of $T$).

\proclaimno Definition. Let $U(T)$ be an affine Hilbert space completion
of\/ $V$ defined above.

There is an obvious $Aut(T)$ action of the group of cellular automorphisms
of $T$ on $U(T)$, with an $Aut(T)$-equivariant isometric embedding of $T$.
An immediate consequence of the construction is

\proclaimno Proposition. If\/ $\Gamma$ acts metrically properly on a tree $T$
then $\Gamma$ is a-T-menable.

In particular $Sl_2\Bbb Z=\Bbb Z/4*_{\Bbb Z/2}\Bbb Z/6$ is a-T-menable.

If $T$ is locally finite,
the (topological) group of all cellular automorphisms
of $T$ is a-T-menable. 
Even if $\Gamma$ acts effectively on $T$, the inclusion 
$\Gamma\to$Aut$(T)$
in general is not closed (therefore we cannot conclude that
$\Gamma$ is a-T-menable). However we have

\proclaimno Proposition. If\/ $\Gamma$ acts on a locally finite tree $T$
and there exists an affine representation of\/ $\Gamma$ on an affine
Hilbert space $W$,
such that stabilizer of any vertex acts properly,
then $U(T)\oplus W$ is $\Gamma$-proper.

\edef\propa{Proposition \tag\ }

{\it Proof :}
Given $g_n\in\Gamma$ with bounded displacement (as acting on $T$),
for any vertex $v$  the distance from $g_nv$ to $v$ is bounded.
Since there are only finitely many such vertices,
one can take a subsequence such that $g_nv$ is constant.
Stabilizers of vertices act properly on $W$, thus
$\{g_n\}$ is relatively compact.
\qed

\proclaimno Proposition. If\/ $\Gamma$ acts on a locally finite tree $T$
and for each vertex $v$ some affine proper representation of $Stab_v$ on
$W_v$ extends to a (perhaps non-proper) affine representation of $\Gamma$,
then $\Gamma$ is a-T-menable.

\edef\propb{Proposition \tag}

{\it Proof :}
This is obvious if the quotient $\Gamma\backslash T$ has finitely
many vertices, since then the sum of the representations corresponding
to any lifts will fulfill the assumptions of \propa

Since there is no infinite sum operation in the category of affine
Hilbert spaces, the construction will depend on choices made.

Let $\{v_k\}$ be a sequence of representatives of the vertices
of the quotients. Let $K_k$ be an exhausting sequence of compact subsets
of $\Gamma$, let $x_k\in W_{v_k}$,
$a_k=k^2+$sup$\{||x_k-gx_k||^2\colon g\in K_k\}$.
Define a norm on $\prod_k W_{v_k}$ in a following way:
$||y-z||^2\colon=\colon\sum a_k^{-2}||y_k-z_k||^2$.
Let $W$ be the completition of the afiine space
$\{y\in \prod_k W_{v_k}\colon ||y-x||<\infty\}$

The choices are made in a such way that $\Gamma$ acts diagonally on $W$
and there are $\Gamma$-equivariant projections (up to scalar change of norm)
to each of $W_{v_k}$.
\qed

{\it Note :} The proof of \propb\ follows the standard
proof of the fact that direct limit of a-T-menable groups is also a-T-menable.

The example of $Sl_2\Bbb Z\sdp\Bbb Z^2$ shows that there are some
obstructions for a representation to extend from the subgroup.
In the terms of group cohomology,
inclusion of groups need not induce epimorphisms
on the level of the first cohomology.

\section Affine representations and subgroups.

If $V$ is $\eu G$-invariant subspace of $W$ then $W/V$ is a linear
$\eu G$-representation (the coset $V$ is a fixed point), therefore
$V$ is metrically proper iff $W$ is.
In general, there is no minimal invariant subspace.

Let $\eu H<\eu G$. 
Assume $V\subset W$ is $\eu H$-invariant subspace.
Fix $x\in W$. Define $\psi_{(V\subset W)}\colon\eu H\backslash G\ni \eu Hg
\mapsto ||gx+V||\in\Bbb R$ (where $gx+V$ is a coset of $gx$ in $W/V$).
If $\psi_{(V\subset W)}$ is proper we say that $(V\subset W)$ is
$\eu H\backslash\eu G$-proper. The definition
does not depend on the choice of $x$.
If $\eu H=\{e\}<\eu G$, and $V$ is any point, then
$(V\subset W)$ is $\eu H\backslash\eu G$-proper exactly
if $W$ is $\eu G$-proper.

\proclaimno Lemma. Let $\eu G_1<\eu G_2<\eu G_3$. If there are
$W_1\subset W_2\subset W_3$, such that $W_i$ is $G_i$ invariant
and $(W_i\subset W_{i+1})$ is $\eu G_i\backslash\eu G_{i+1}$-proper
then $(W_1\subset W_3)$ is $\eu G_1\backslash\eu G_3$-proper.

\edef\lema{Lemma \tag\ }

{\it Proof :} If $\psi_{(W_1\subset W_3)}(\eu G_1 g_n)$ is bounded,
then $\psi_{(W_2\subset W_3)}(\eu G_2 g_n)$ is bounded.
Therefore $g_n=g'_nh_n$ where $g'_n\in\eu G_2$ and $\{h_n\}$
is relatively compact. By the triangle inequality $||g'_nx-x||\leq
||g_n(h^{-1}_nx-x)||+||g_nx-x||+||h^{-1}_nx-x||$,
so $\psi_{(W_1\subset W_2)}(\eu G_1g'_n)$ is bounded,
therefore $g'_n=g''_nh'_n$,
where $g''_n\in\eu G_1$ and $\{h'_n\}$ is relatively compact.
Finally $g_n=g''_n(h'_nh_n)$.\qed

Unfortunately it is not known whether for any $\eu H<\eu G$
and proper $\eu G$-representation $W$ there exist $\eu H$-invariant
subspace $V$, such that  $(V\subset W)$ is
$\eu H\backslash\eu G$-proper.
The cases when it does happend are discussed
in the following sections.

\section Free products with amalgamation. 

If $\Gamma=\eu G_1*_{\eu H}\eu G_2$
then, according to Serre theory [S], a graph $T$ with the
set of vertices equal to
$\Gamma/\eu G_1\cup\Gamma/\eu G_2$ and the set
of edges equal to $\Gamma/\eu H$, with inclusion as incidence relation,
is a tree (on which $\Gamma$ acts on the left).
The representations of $\eu G_1$ and $\eu G_2$ on $U(T)$ have global
fixed points.

Let $\eu H$ be a common subgroup in $\eu G_1$ and $\eu G_2$.
Let $W_i$ be $\eu G_i$-representations.
Let $W$ be their common $\eu H$-invariant subspace.
Define $\Cal H_i\colon=W_i/W$.
Inductively decompose
$\uparrow^{\eu G_i}_{\eu H}\Cal H_{\omega}$
(where $\omega$ is a sequence of 1s and 2s)
with respect to $\eu H$ as $\Cal H_\omega\oplus\Cal H_{i\omega}$.
$\eu G_1$ acts on
$\Cal H_2^\bullet=(\Cal H_{2}\oplus\Cal H_{12})\oplus(\Cal H_{212}\oplus\Cal H_{1212})\oplus\dots$
and $\Cal H_1^\circ=(\Cal H_{21}\oplus\Cal H_{121})\oplus\dots$,
$\eu G_2$ acts on
$\Cal H_1^\bullet=(\Cal H_{1}\oplus\Cal H_{21})\oplus\cdots$
and $\Cal H_2^\circ=(\Cal H_{12}\oplus\Cal H_{212})\oplus\cdots$.
Both representations of $\eu H$ on $\Cal H_i^\bullet$ coincide.

\proclaimno Definition.
Let $W_\Gamma=W\oplus\Cal H_1^\bullet\oplus\Cal H_2^\bullet
=W_1\oplus\Cal H_1^\circ\oplus\Cal H_2^\bullet
=W_2\oplus\Cal H_1^\bullet\oplus\Cal H_2^\circ$.

An immediate consequence from the construction is

\proclaimno Theorem. Let $\eu H<\eu G_i$, $i=1,2$.
Let $W_i$ be $\eu G_i$-representations.
Let $W$ be their common $\eu H$-invariant subspace.
Let $\Gamma=\eu G_1*_{\eu H}\eu G_2$. Then 
$W_i$ and $W$ are respectively $\eu G_i-$ and $\eu H-$invariant 
subspaces of $W_\Gamma$.

\edef\thma{Theorem \tag\ }

{\it Note :} Although there is no way to induce an affine
representation, $W_\Gamma$ is morally equal to
$\uparrow_{\eu G_1}^\Gamma W_1\oplus\uparrow_{\eu G_2}^\Gamma W_2/
\uparrow_{\eu H}^\Gamma W$. 
If $W'\subset W$ is another $\eu H$-invariant subspace, then $W_\Gamma$
is not a $\Gamma$-invariant subspace of $W'_\Gamma$ (constructed from
the triple $W_1\supset W'\subset W_2$).  

A straightforward consequence of \propa and \thma is the following

\proclaimno Corollary. If\/ $\eu H$ is of finite index in $\eu G_1$
and $\eu G_2$ and if there are metrically proper representations
of\/  $\eu G_i$ that coincide when restricted to $\eu H$, then
$\eu G_1*_{\eu H}\eu G_2$ is a-T-menable.

{\it Example [JJV,BCS]. } The torus group
$\Gamma_{p,q}=\langle x,y\vert x^p=y^q\rangle$ is a-T-menable. 

\proclaimno Theorem. With the assumptions of \thma,
if\/ $W_i$ is proper affine $\eu G_i$-represen\-tation and
$(W\subset W_i)$ is $\eu H\backslash \eu G_i$-proper for $i=1,2$
then\hfill\break
(1) $W_{\Gamma}\oplus U(T)$ is $\Gamma$-proper,\hfill\break
(2) if\/ $\eu H'<\eu G_1$ and $V$ is proper $\eu H'$-space such that
$(V\subset W_1)$ is $\eu H'\backslash \eu G_1$-proper, then
$(V\subset W_{\Gamma}\oplus U(T))$ is $\eu H'\backslash\Gamma$-proper.

\edef\thmb{Theorem \tag\ }

{\it Proof :}
In fact (1) is a special case of (2). Therefore we will prove (2).

Given $\gamma_n\in\Gamma$ such that
$\psi_{(W,V)}(\gamma_n)$ is bounded.
Define {\sl the length function} $\ell\colon\Gamma\to\Bbb N$
by $\ell_{|\eu H}\equiv 0$,
$\ell(\gamma)=min\{\ell(\eta)+1\vert\gamma=\eta g$, where 
$g\in \eu G_1\cup \eu G_2\}$. This function is equal to the distortion
of the action on $U(T)$. Therefore we may find a subsequence such that
$\ell(\gamma_n)=k$.

If $k\leq 1$ there is nothing to prove. If all but finitely many
$\gamma_n\in\eu G_2$ we have to use \lema

Let $x\in V$.
Define $\varphi_x(\gamma)$ to be the component of $\gamma x$ in
$\bigoplus_{|\omega|=l(\gamma)}\Cal H_\omega$.

{\advance\leftskip by \bigskipamount\advance\rightskip by \bigskipamount

\proclaimno Lemma. If $\gamma=\gamma_1\gamma_2$,
$\ell(\gamma)=\ell(\gamma_1)+\ell(\gamma_2)$ and $\ell(\gamma_2)\geq1$,
then $||\varphi_x(\gamma)||=||\varphi_x(\gamma_2)||$.

\global\edef\lemb{Lemma \tag\ }

{\it Proof :} Without loss of generality $\ell(\gamma_1)=1$.
From the definition of induced representation
$\gamma_1\Cal H_\omega\perp\Cal H_\omega$, therefore
$\varphi_x(\gamma)=\gamma_1\varphi_x(\gamma_2)$.\qed

}

Now we proceed by induction on $k$ as follows:
we define $\eta_n$ such that
$\gamma_n=\eta_ng_n$ ($g_n\in \eu G_i$) and $l(\eta_n)=k-1$.

From \lemb we see $||\gamma_nx-x||\geq
||\varphi_x(\gamma_n)||=||\varphi_x(g_n)||=\psi_{(W_i,W)}(g_n)$.
Therefore $\{g_n\}$ is relatively compact.
Since $\eta_nx=\gamma_nx-\eta_n(g_nx-x)$
and, by induction assumption, $\psi_{(V,W_\Gamma)}$ is proper
when restricted to cosets of elements of length smaller than $k$,
we obtain the claim.\qed

\section Groups that act on trees with finite edge stabilizers.

The first case, where it is easy to fulfill the assumptions of the
\thmb is when $\eu H$ is finite, since then one can
find a fixed point of any $\eu H$-representation
simply taking the center of mass $*$ of any orbit.
The pair $(\{*\}\subset W_i)$ is $\eu H\backslash \eu G_i$ proper
iff $W_i$ is proper $\eu G_i$-space. Summarizing this:

\proclaim Proposition. [JJV, Prop.~6.2.3 (1)] Let $\eu G_1$, $\eu G_2$ be two
grups containing finite subgroup $H$,
and let $\Gamma=\eu G_1*_{\eu H}\eu G_2$ be the corresponging amalgamated
product. If $\eu G_1$ and $\eu G_2$ are a-T-menable, then so is $\Gamma$.

\jjv is an easy consequence of the above [JJV].

\section Baumslag-Solitar groups ant their certain generalizations

The second easy case occurs when $\eu H$ is of finite index in $\eu G_i$'s.
Then any pair is automatically $\eu H\backslash \eu G_i$ proper.

Let us recall some definitions from [GJ].
Let $\eu G\subset\frakN$ be a closed subgroup of
a locally compact compactly generated topological group $\frakN$.
Let $i_k\colon\eu H\to\eu G$ ($k=1,2$)
be two inclusions onto finite index open subgroups, which are conjugated
by an automorphism $\phi$ of $\frakN$.

\proclaimno Definition. The $\frakN$-BS group $\Gamma$ is the
group derived from $(\eu G,\eu H,i_1,i_2)$ by the
(topological) HNN construction.

\proclaimno Theorem. [GJ] If\/ $\frakN$ is a-T-menable then $\frakN$-BS groups
are a-T-menable.

{\it Proof :}
We mimic the proof [JJV 6.2.7] for the case of a
HNN extension, where the edge stabilizer is finite.

{\it Step 1.} Let $\Gamma_0$ be a fundamental group of the following tree:
$$\matrix{&\eu G&&&&\eu G&&&&\eu G&\cr
\cdots&&\nwarrow\rlap{$i_1$}&&\nearrow\rlap{$i_2$}&&\nwarrow\rlap{$i_1$}&&\nearrow\rlap{$i_2$}&&\cdots\cr
&&&\eu H&&&&\eu H&&&\cr}$$

We have to find consistent representations of different copies of $\eu G$.
The Hilbert space in each case will be the one on which $\frak N$
acts
properly. The $k$-th copy acts by $\phi^k(\cdot)$, where $\phi$ is
the automorphism of $\frak N$ that conjugates $i_1$ and $i_2$.

By induction, each of
$$\matrix{\eu G*_{\eu H}\cdots*_{\eu H}\eu G&&&&\eu G\cr
&\nwarrow\rlap{$i_1$}&&\nearrow\rlap{$i_2$}&\cr
&&\eu H&&\cr}$$
satisfies assumptions of \thmb (alternatively:
by \thma and induction we construct representation of $\Gamma_0$
and then use \propb).
It is easy to show [JJV Prop.~6.1.1] that an increasing union of
open a-T-menable subgroups is again a-T-menable.
Therefore $\Gamma_0$ is a-T-menable.

{\it Step 2.} $\Gamma=\Gamma_0\sdp\Bbb Z$, where $\Bbb Z$ acts
through the shift. Therefore $\Gamma$ is an extension of an a-T-menable group
with amenable quotient, so $\Gamma$ is a-T-menable by \joli.\qed

\section References.

\item{[BCS]} C.~B\'eguin and T. Ceccherini-Silberstein, {\it Formes faibles
de moyenabilit\'e pour les gruupes \`a un relateur},
Bull.~Belg.~Math.~Soc.~Simon Stevin, {\bf 1} (2000), pp.~135-148

\item{[GJ]} S.~R.~Gal, T.~Januszkiewicz, {\it New a-T-menable HNN-extensions},
J. Lie Theory Vol. {\bf 13} (2003), No. 2, pp.~383--385

\item{[H]} U.~Haagerup, {\it An example of non-nuclear $C^*$-algebra which has
the metric approximation property}, Invent.~Math. {\bf 50} (1979), pp.~279-293

\item{[HV]} P.~de la Harpe, A.~Valette, {\it La properi\'et\'e (T) de Kazhdan
pour les groupes localement compacts}, Asterisque {\bf 157}, 1989

\item{[JJV]} P.~Jolissant, P.~Julg, A.~Valette,
{\it Discrete groups}, Chapter 6.~from
{\it Groups with the Haagerup property (Gromov's a-T-menability)},
Birkh\"auser Verlag, 2001

\item{[S]} J.~P.~Serre, {\it Trees}, Springer-Verlag, 1980

\bigskip

\rightline{\it Wroc\l aw, December 2000}

\bye